  \documentclass[smallcondensed]{svjour3}

  \usepackage[letterpaper, left=1in, right=1in, top=1in, bottom=1in]{geometry}
  
  \usepackage[affil-it]{authblk}  
  
  \usepackage{palatino}  
  \usepackage{setspace}  
  \usepackage{color}  
  \usepackage[dvipsnames]{xcolor}  

	\usepackage[labelfont=bf]{caption}	

  \usepackage{amsmath}  
  \usepackage{amssymb}  
  \usepackage{bbm}  
  \allowdisplaybreaks[4]  
	\usepackage{relsize}  
	\usepackage{algorithm}  
	\usepackage{algpseudocode}  
  
  \DeclareMathOperator*{\minimize}{minimize}
  
  \newcommand{\st}{\mathrm{subject\;to}}

  \usepackage{pifont}  

	\usepackage{longtable}  
	\usepackage{multirow}  
	\usepackage{array}  
	\usepackage{booktabs}  
    
  \usepackage{graphicx}  
  \graphicspath{{Figures/}}  
  \usepackage{float}  
  \usepackage{subfig}  
	\usepackage{hyperref}
	\hypersetup{
		pdfstartview = {FitV},  
		colorlinks = true,    
		linkcolor = NavyBlue,  
		citecolor = NavyBlue,  
		filecolor = NavyBlue,  
		urlcolor = NavyBlue  
	}  

  \usepackage[square,numbers]{natbib}  
\title{Branch-and-Price for a Class of Nonconvex Mixed-Integer Nonlinear Programs}
\author{Andrew Allman \and Qi Zhang}
\institute{Andrew Allman \and Qi Zhang (corresponding author) \at 
	Department of Chemical Engineering and Materials Science\\ University of Minnesota, Twin Cities\\ Minneapolis, MN 55455, USA
	\\Email: qizh@umn.edu}
\date{}
\titlerunning{Branch-and-Price for MINLPs}
\authorrunning{Allman and Zhang}
\journalname{Journal of Global Optimization}

\begin{document}

\maketitle

\begin{abstract}
This work attempts to combine the strengths of two major technologies that have matured over the last three decades: global mixed-integer nonlinear optimization and branch-and-price. We consider a class of generally nonconvex mixed-integer nonlinear programs (MINLPs) with linear complicating constraints and integer linking variables. If the complicating constraints are removed, the problem becomes easy to solve, e.g. due to decomposable structure. Integrality of the linking variables allows us to apply a discretization approach to derive a Dantzig-Wolfe reformulation and solve the problem to global optimality using branch-and-price. It is a remarkably simple idea; but to our surprise, it has barely found any application in the literature. In this work, we show that many relevant problems directly fall or can be reformulated into this class of MINLPs. We present the branch-and-price algorithm and demonstrate its effectiveness (and sometimes ineffectiveness) in an extensive computational study considering multiple large-scale problems of practical relevance, showing that, in many cases, orders-of-magnitude reductions in solution time can be achieved.
\end{abstract}

\section{Introduction}

Mixed-integer nonlinear programming (MINLP) has proven to be a powerful modeling paradigm and has received increased attention in recent years. However, despite the tremendous advances in the theoretical and algorithmic treatment of MINLPs, they are still significantly less scalable than their linear counterparts such that solving MINLPs of large sizes remains a challenge. In this work, we consider the exact solution of a class of generally nonconvex MINLPs whose structure is amenable to branch-and-price, hence allowing us to combine the strengths of two key enabling technologies: column generation for large-scale integer programming and global optimization of MINLPs. Specifically, we consider problems of the following form:
\begin{subequations}
\label{eqn:MINLP}
\begin{align}
\minimize_{x,y,z} \quad & c^{\top} x + f(y,z) \label{eqn:MNILPobj} \\
\st \quad & A x + D y \geq b \label{eqn:MINLPlinCon} \\
& g(y,z) \leq 0 \label{eqn:MINLPnlinCon} \\
& y^{\min} \leq y \leq y^{\max} \label{eqn:MINLPbounds} \\
& x \in \mathbb{R}_+^m \times \mathbb{Z}_+^{\bar{m}}, \, y \in \mathbb{Z}^p, \, z \in \mathbb{R}^q \times \mathbb{Z}^{\bar{q}},
\end{align}
\end{subequations}
where the vectors $x$ and $z$ can contain both continuous and integer variables, while the $y$-variables are all integer. The objective function \eqref{eqn:MNILPobj} consists of a linear term, $c^{\top}x$, and a nonlinear term, $f(y,z)$. While the linear constraints \eqref{eqn:MINLPlinCon} involve $x$ and $y$, $g(y,z)$ in \eqref{eqn:MINLPnlinCon} are nonlinear functions of $y$ and $z$. The functions $f$ and $g$ can be nonconvex. Constraints \eqref{eqn:MINLPbounds} ensure that $y$ are bounded. We assume that problem \eqref{eqn:MINLP} can be efficiently solved if constraints \eqref{eqn:MINLPlinCon}, which we call the \textit{complicating constraints}, are removed.

We are particularly interested in problems that, without the complicating constraints, decompose into smaller independent subproblems that can be solved using state-of-the-art global MINLP solvers. We find that many relevant problems directly fall or can be reformulated into this class of MINLPs. Examples can be found in process and product design, production capacity planning, dynamic facility location, stochastic programming, statistical learning, and many more application domains.

In this work, we investigate the computational feasibility of a branch-and-price approach to solving MINLPs of the given form. While the suitability of the proposed algorithm has been indicated in the literature, it has barely found any application. We hence aim to systematically analyze the theoretical basis of the branch-and-price approach, highlight critical algorithmic considerations, and examine the algorithm's performance in computational experiments.

The remainder of this paper is organized as follows. First, related literature is reviewed in Section \ref{sec:Literature}. Next, we present a discretization approach for reformulating problems of the form \eqref{eqn:MINLP} in Section \ref{sec:Reformulation}. This reformulation makes the problem amenable for solution using branch-and-price, and the algorithm for doing so is presented in Section \ref{sec:Algorithm}. We show how this approach can be extended to problems where the pricing subproblem is decomposable in Section \ref{sec:Decompose}. In Section \ref{sec:Experiments}, the computational performance of the proposed branch-and-price algorithm is compared to solving the full-space problem using a state-of-the-art global MINLP solver in four representative case studies. Finally, in Section \ref{sec:Conclusions}, we provide some concluding remarks.

\section{Literature Review}
\label{sec:Literature}

The development of exact methods for the solution of mixed-integer \textit{linear} programs (MILPs) dates back to the 1950s \cite{Dantzig1954, Gomory1958} (for more details on the history of integer programming, see \cite{Junger2009}). Over the last decades, MILP has reached a level of maturity that has made it the primary approach to solving many industrial and scientific problems of high complexity and dimensionality. Mixed-integer \textit{nonlinear} programming is a more recent development and was initially motivated by applications in chemical and process systems engineering \cite{Grossmann1979}.

The development of MINLP algorithms has initially focused on the convex case, i.e. problems in which, loosely speaking, the continuous relaxation of the MINLP is convex. Methods for solving convex MINLPs include branch-and-bound \cite{Gupta1985, Stubbs1999}, generalized Benders decomposition \cite{Geoffrion1972}, outer approximation \cite{Duran1986, Fletcher1994}, LP/NLP-based branch-and-bound \cite{Quesada1992}, extended cutting plane \cite{Westerlund1995}, and extended supporting hyperplane \cite{Kronqvist2016}. For recent reviews on convex MINLP, we refer the reader to \cite{Grossmann2002}, \cite{Bonami2012}, and \cite{Kronqvist2019}.

Compared to convex MINLPs, solving nonconvex MINLPs is significantly more challenging due to the nonconvexity that remains even after relaxing the integer restrictions. Exact algorithms for nonconvex MINLP incorporate concepts from global continuous optimization, such as convex relaxations and spatial branch-and-bound \cite{McCormick1976}. Major improvements have been achieved with the development of the branch-and-reduce \cite{Ryoo1996, Tawarmalani2004} and $\alpha$-branch-and-bound methods \cite{Androulakis1995, Adjiman2000}. Furthermore, incorporating MILP relaxations and integer programming techniques for generating cutting planes has proven to be very effective \cite{Kilinc2018}. These and other algorithmic advances are implemented in state-of-the-art global MINLP solvers, such as BARON \cite{Kilinc2018}, Couenne \cite{Belotti2009}, LINDOGlobal \cite{Lin2009}, ANTIGONE \cite{Misener2014}, and SCIP \cite{Vigerske2018}. For recent reviews focusing on nonconvex MINLP, see \cite{Burer2012} and \cite{Boukouvala2016}.

Although the performance of global MINLP solvers has improved significantly over the last two decades, they are still by far not as scalable as state-of-the-art MILP solvers. Hence, to solve large-scale nonconvex MINLPs, one often resorts to decomposition methods that exploit special model structures. Popular decomposition approaches include different variants of Lagrangean decomposition \cite{Guignard2003} and progressive hedging \cite{Watson2011}. However, these methods have to be considered heuristics since, although they can provide lower and upper bounds, it is not guaranteed that the duality gap can be closed. Exact decomposition algorithms for nonconvex MINLP are a rarity. One of those is bilevel decomposition, which iterates between a master MILP that is a relaxation of the original MINLP and an NLP or MINLP subproblem obtained by fixing integer variables. Integer and outer-approximation cuts (or tailored cuts that can be interpreted as such) are added to the master MILP at each iteration. The convergence behavior of bilevel decomposition strongly depends on the quality of the MILP relaxation and is therefore very application-specific \cite{Lotero2016, Lara2018a, Elsido2019}. Generalized Benders decomposition has been extended to solve two-stage stochastic nonconvex separable MINLPs in which only the continuous variables are involved in the nonconvex terms \cite{Li2011b}. A branch-and-bound scheme to address two-stage stochastic nonconvex MINLP with mixed-integer first-stage and continuous second-stage variables is proposed in \cite{Cao2019}. Combining generalized Benders decomposition and branch-and-cut, \cite{Li2019a} considers the case with nonconvex constraints and mixed-binary variables in both stages. Finally, \cite{Rebennack2009} presents a decomposition method based on column enumeration for nonconvex MINLPs with an assignment constraint; however, since the number of columns grows exponentially with the number of assignment decisions, this method is only suited for problems with a few assignment variables.

In the same spirit of decomposition, we apply in this work branch-and-price, which has its origin in the pioneering works of \cite{Dantzig1960}, who introduced the fundamental idea of column generation for linear programming, and \cite{Desrosiers1984}, who were the first to embed column generation in a branch-and-bound framework to solve a large-scale MILP. Branch-and-price has been a success story in large-scale mixed-integer optimization, with applications in vehicle routing \cite{Desrochers1992, Desaulniers2002}, crew scheduling \cite{Desrochers1989, Stojkovic1998}, fleet assignment \cite{Ioachim1999, Belanger2006}, and capacity planning \cite{Flores-Quiroz2019}, to just name a few. For reviews on branch-and-price, see \cite{Barnhart1998} and \cite{Lubbecke2005}.

The vast majority of existing works on branch-and-price consider integer and mixed-integer linear problems. However, there is significant flexibility in incorporating nonlinearities in the pricing problem, which has been mentioned (in form of brief side notes) in the literature \cite{Lubbecke2005, Singh2009}. It is all the more surprising that there seems to be almost no published work on applying branch-and-price with mixed-integer nonlinear pricing problems. The only notable exception that we have been able to find is \cite{Nowak2018}, which proposes a column-generation-based method for generating inner and outer approximations for nonconvex MINLPs. We believe that there is considerable room for theory and algorithm development in the application of branch-and-price concepts to solving MINLPs.


\section{Reformulation via Discretization}
\label{sec:Reformulation}

\begin{figure}[ht]\centering
\includegraphics[width=3in]{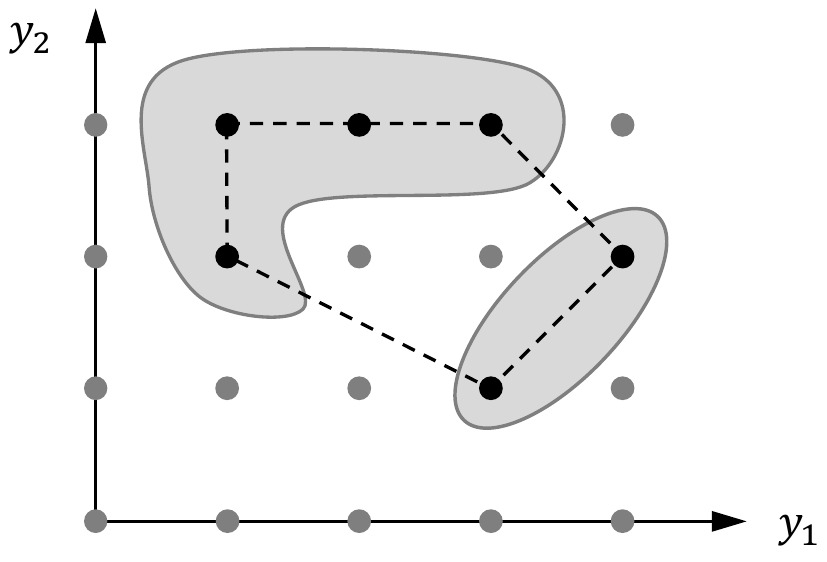}
\caption{The two-dimensional feasible set $\mathcal{Y}$ is given by the set of dark-colored points. The shaded area represents set $\mathcal{Y}$ without the integrality restrictions on $y$.}
\label{fig:Discretization}
\end{figure}

To make the problem amenable to branch-and-price, we first apply a discretization approach \cite{Vanderbeck2000} to derive an extensive formulation of problem \eqref{eqn:MINLP}. Consider the feasible set for $y$ when disregarding constraints \eqref{eqn:MINLPlinCon}:
\begin{equation}
\mathcal{Y} := \left\lbrace y: g(y,z) \leq 0, \; y^{\min} \leq y \leq y^{\max}, \; y \in \mathbb{Z}^p, \; z \in \mathbb{R}^q \times \mathbb{Z}^{\bar{q}} \right\rbrace,
\end{equation}
which is a finite set of integer points with cardinality $K := |\mathcal{Y}|$. Hence, $\mathcal{Y}$ can be equivalently expressed as:
\begin{equation}
\mathcal{Y} = \{\bar{y}_1, \dots, \bar{y}_K\} = \left\lbrace y: y = \sum_{k \in \mathcal{K}} \lambda_k \, \bar{y}_k, \; \sum_{k \in \mathcal{K}} \lambda_k = 1, \; \lambda_k \in \{0,1\} \; \forall \, k \in \mathcal{K} \right\rbrace,
\end{equation}
where $\bar{y}_k$ denotes a specific point in $\mathcal{Y}$ and $\mathcal{K} := \{1,\dots,K\}$. Here, the binary variable $\lambda_k$ is equal to 1 if $y = \bar{y}_k$ and otherwise 0. Figure \ref{fig:Discretization} shows a two-dimensional example of $\mathcal{Y}$, which is given by the dark-colored points. The shaded area, which is clearly nonconvex, represents $\mathcal{Y}$ without the integrality restrictions on $y$.

\begin{remark}
Note that in the example shown in Figure \ref{fig:Discretization}, $\mathcal{Y} \neq \mathrm{conv}(\mathcal{Y}) \cap \mathbb{Z}^p$, where $\mathrm{conv}(\mathcal{Y})$ denotes the convex hull of $\mathcal{Y}$. In general, $\mathcal{Y} \subseteq \mathrm{conv}(\mathcal{Y}) \cap \mathbb{Z}^p$, which is due to nonconvex constraint functions $g$ or integer components in $z$. Hence, the traditional convexification technique \cite{Lubbecke2005} that is often used in Dantzig-Wolfe reformulations does not apply here.
\end{remark}

For each $k \in \mathcal{K}$, we can define an optimal cost:
\begin{equation}
\bar{f}_k := \min_{z \in \mathbb{R}^q \times \mathbb{Z}^{\bar{q}}} \left\lbrace f(\bar{y}_k,z): g(\bar{y}_k,z) \leq 0 \right\rbrace,
\end{equation}
which is the last ingredient that we need to reformulate \eqref{eqn:MINLP} into the following \textit{master problem}:
\begin{subequations}
	\label{eqn:MP}
\begin{align}
\mathrm{(MP):} \quad v^{\mathrm{MP}} := \min_{x,\lambda} \;\; & c^{\top} x + \sum_{k \in \mathcal{K}} \lambda_k \bar{f}_k \\
\mathrm{s.t.} \;\; & A x + D \sum_{k \in \mathcal{K}} \lambda_k \, \bar{y}_k \geq b \label{eqn:MPCon} \\
& \sum_{k \in \mathcal{K}} \lambda_k = 1 \label{eqn:ConvexityCon} \\
& x \in \mathbb{R}_+^m \times \mathbb{Z}_+^{\bar{m}} \\
& \lambda_k \in \{0,1\} \quad \forall \, k \in \mathcal{K},
\end{align}
\end{subequations}
which is a large-scale MILP as $K$ grows exponentially with the dimension of $y$. Problems (MP) and \eqref{eqn:MINLP} are equivalent in a sense that they have the same optimal value, $v^{\mathrm{MP}}$, and every solution of (MP) can be mapped to a unique solution of \eqref{eqn:MINLP} in the $(x,y)$-space and a corresponding set of $z$-values that provide the same optimal value (the converse is trivially true).

\section{The Branch-and-Price Algorithm}
\label{sec:Algorithm}

In the branch-and-price algorithm, we solve (MP) using branch-and-bound, but solve the LP relaxation at each node via column generation. In the following, we describe the major elements of the algorithm.

\subsection{Column Generation}

Consider the LP relaxation of (MP), which we denote by $(\overline{\mathrm{MP}})$. The large number of variables in $(\overline{\mathrm{MP}})$ typically prohibits solving it in full space. Instead, we consider a \textit{restricted master problem}, denoted by (RMP), which involves only a subset of columns $\widehat{\mathcal{K}} \subseteq \mathcal{K}$. The LP relaxation of (RMP), denoted by $(\overline{\mathrm{RMP}})$, can be solved efficiently as long as the size of $\widehat{\mathcal{K}}$ is sufficiently small. New columns are generated as needed by solving the following \textit{pricing problem}:
\begin{subequations}
\begin{align}
\mathrm{(PP):} \quad \zeta := \min_{y,z} \;\; & f(y,z) - \pi^{\top} D y - \mu \\
\mathrm{s.t.} \;\; & g(y,z) \leq 0 \\
& y^{\min} \leq y \leq y^{\max} \\
& y \in \mathbb{Z}^p, \, z \in \mathbb{R}^q \times \mathbb{Z}^{\bar{q}},
\end{align}
\end{subequations}
where $\pi$ and $\mu$ denote the values of the dual variables associated with constraints \eqref{eqn:MPCon} and \eqref{eqn:ConvexityCon}, respectively, at the optimal solution of $(\overline{\mathrm{RMP}})$. Problem (PP) minimizes the reduced cost; hence, if $\zeta < 0$, the corresponding $y^*$ may improve the solution to $(\overline{\mathrm{MP}})$ and is therefore added as a new column to $\widehat{\mathcal{K}}$. Note that (PP) is a generally nonconvex MINLP.

Since $(\overline{\mathrm{RMP}})$ considers a restricted set of columns, its optimal value, $v^{\overline{\mathrm{RMP}}}$, is an upper bound on the optimal value of $(\overline{\mathrm{MP}})$, $v^{\overline{\mathrm{MP}}}$. Following standard duality arguments (\cite{Wolsey1998}, p. 189), one can show that a lower bound is given by $v^{\overline{\mathrm{RMP}}} + \zeta$. Hence, we have the following relationship:
\begin{equation}
v^{\overline{\mathrm{RMP}}} + \zeta \leq v^{\overline{\mathrm{MP}}} \leq v^{\overline{\mathrm{RMP}}}.
\end{equation}

Algorithm \ref{alg:solveRelaxedMP} shows the pseudocode of the algorithm for solving $(\overline{\mathrm{MP}})$. Here, the current lower and upper bounds on $v^{\overline{\mathrm{MP}}}$ are denoted by $\mathrm{LB}^{\overline{\mathrm{MP}}}$ and $\mathrm{UB}^{\overline{\mathrm{MP}}}$, respectively. The column generation algorithm is finite and exact. If we set the tolerance $\bar{\epsilon}$ to zero, $\zeta$ will be zero at the last iteration and the algorithm terminates at the optimal solution of $(\overline{\mathrm{MP}})$.

\begin{algorithm}[h]
\caption{Column generation algorithm for solving $(\overline{\mathrm{MP}})$.}
\begin{algorithmic}[1]
\Function{solveRelaxedMP}{$(\overline{\mathrm{MP}})$, (PP), $\widehat{\mathcal{K}}$}
  \State Set tolerance $\bar{\epsilon}$
  \State Initialize: $\mathrm{LB}^{\overline{\mathrm{MP}}} \gets -\infty$, $\mathrm{UB}^{\overline{\mathrm{MP}}} \gets \infty$, $k \gets |\widehat{\mathcal{K}}|$
  \While{$\mathrm{UB}^{\overline{\mathrm{MP}}} - \mathrm{LB}^{\overline{\mathrm{MP}}} > \bar{\epsilon}$}
    \State Solve $(\overline{\mathrm{RMP}})$, obtain $x^*$, $\lambda^*$, $\pi$, $\mu$, and $v^{\overline{\mathrm{RMP}}}$
    \State Update upper bound: $\mathrm{UB}^{\overline{\mathrm{MP}}} \gets v^{\overline{\mathrm{RMP}}}$
    \State Solve (PP), obtain $y^*$ and $\zeta$
    \If{$\zeta < 0$} add column
      \State $k \gets k+1$
      \State $\bar{y}_k \gets y^*$, $\bar{f}_k \gets \zeta + \pi^{\top} D y^* + \mu$, and $\widehat{\mathcal{K}} \gets \widehat{\mathcal{K}} \cup \{k\}$
    \EndIf
    \State Update lower bound: $\mathrm{LB}^{\overline{\mathrm{MP}}} \gets \max \{\mathrm{LB}^{\overline{\mathrm{MP}}}, v^{\overline{\mathrm{RMP}}} + \zeta\}$
	\EndWhile
  \State \Return{$x^*$, $\lambda^*$, $v^{\overline{\mathrm{MP}}} \gets v^{\overline{\mathrm{RMP}}}$}
\EndFunction
\end{algorithmic}
\label{alg:solveRelaxedMP}
\end{algorithm}

\begin{remark}
\label{rmk:SuboptimalPP}
If we solve (PP) using a global MINLP solver such as BARON, we obtain, if available, lower ($l$) and upper ($u$) bounds on $\zeta$ at every iteration of the branch-and-bound algorithm. Using these bounds, we can generate new columns and compute valid bounds on $v^{\overline{\mathrm{MP}}}$ without solving (PP) to optimality as follows: A column is added if $u < 0$. In line 10 of Algorithm \ref{alg:solveRelaxedMP}, $\bar{f}_k$ is computed using $u$, i.e. $\bar{f}_k \gets u + \pi^{\top} D y^* + \mu$. Then the following bounds can be computed:
\begin{equation}
v^{\overline{\mathrm{RMP}}} + l \leq v^{\overline{\mathrm{RMP}}} + \zeta \leq v^{\overline{\mathrm{MP}}} \leq v^{\overline{\mathrm{RMP}}}.
\end{equation}
Only in the last iteration, (PP) has to be solved to $\bar{\epsilon}$-optimality in order to achieve a desired optimality gap. Not solving (PP) to optimality at every iteration can mitigate the impact of the tailing-off effect in solving the MINLP and significantly speed up the overall algorithm.
\end{remark}

\subsection{Branch-and-Bound}

If the solution of $(\overline{\mathrm{MP}})$ is integer feasible, it is also the optimal solution to (MP); remarkably, this is often the case, as our computational experiments show (see Section \ref{sec:Experiments}). In general, however, the solution may not satisfy the integrality restrictions, in which case we have to apply branch-and-bound. It is well known that branching on the $\lambda$-variables leads to unbalanced branch-and-bound trees; hence, it is recommended to design branching rules based on the original variables \cite{Vanderbeck2000}. We comment, whenever appropriate, on branching rules in the case studies in Section \ref{sec:Experiments} as they often have to be tailored to the specific application.

An outline of the branch-and-price algorithm is shown in Algorithm \ref{alg:BnP}. Here, $\mathcal{N}$ denotes the set of nodes generated in the branch-and-bound tree, which is initialized with the root node 1. Furthermore, we start with a nonempty set of feasible columns $\widehat{\mathcal{K}}$. At each node $n$, we solve the linear relaxation of the master problem associated with that node, denoted by $(\overline{\mathrm{MP}})_n$, using the corresponding pricing problem (PP)$_n$. The solution of $(\overline{\mathrm{MP}})_n$ can be used to update the overall lower bound $\mathrm{LB}^{\mathrm{MP}}$, and if it is integer feasible, it also provides an upper bound on $v^{\mathrm{MP}}$; otherwise, we can solve (RMP)$_n$, which considers the integrality constraints, to obtain a feasible solution to (MP) and hence an upper bound. The incumbent solution $(\hat{x}^*, \hat{\lambda}^*)$ is updated every time a new feasible solution is found. The algorithm terminates when an optimality gap smaller than $\epsilon$ is reached; otherwise, branching rules are applied to update the node set $\mathcal{N}$, and the next node $n$ to be evaluated is selected.

\begin{algorithm}[h]
\caption{Branch-and-price algorithm for solving $(\mathrm{MP})$.}
\begin{algorithmic}[1]
\State Set tolerance $\epsilon$
\State Initialize: $\mathrm{LB}^{\mathrm{MP}} \gets -\infty$, $\mathrm{UB}^{\mathrm{MP}} \gets \infty$, $\mathcal{N} \gets \{1\}$, $n \gets 1$, $\widehat{\mathcal{K}}$
\While{$|\mathcal{N}| > 0$}
  \State \Call{solveRelaxedMP}{$(\overline{\mathrm{MP}})_n$, (PP)$_n$, $\widehat{\mathcal{K}}$}
  \State Update lower bound $\mathrm{LB}^{\mathrm{MP}}$
  \If{solution not integer feasible}
    \State Solve (RMP)$_n$, obtain $x^*$, $\lambda^*$, $v^{\mathrm{RMP}}$
  \EndIf
  \State Update upper bound $\mathrm{UB}^{\mathrm{MP}}$ and incumbent solution $(\hat{x}^*, \hat{\lambda}^*)$
  \If{$\mathrm{UB}^{\mathrm{MP}} - \mathrm{LB}^{\mathrm{MP}} > \epsilon$}
    \State Apply branching rules and update $\mathcal{N}$
    \State Select next node $n$
  \Else
    \State \Return{$\hat{x}^*$, $\hat{\lambda}^*$, $\mathrm{LB}^{\mathrm{MP}}$, $\mathrm{UB}^{\mathrm{MP}}$}
  \EndIf
\EndWhile
\end{algorithmic}
\label{alg:BnP}
\end{algorithm}

\begin{remark}
As mentioned in Remark \ref{rmk:SuboptimalPP}, a lower bound on $v^{(\overline{\mathrm{MP}})_n}$ is obtained at every iteration of \Call{solveRelaxedMP}{$(\overline{\mathrm{MP}})_n$, (PP)$_n$, $\widehat{\mathcal{K}}$}. This information can be used to potentially prune the node before the column generation algorithm terminates.
\end{remark}

\begin{remark}
Solving (RMP)$_n$ after solving $(\overline{\mathrm{MP}})_n$ at every node (line 7 of Algorithm \ref{alg:BnP}) is not required for convergence; however, it provides upper bounds that can help reduce the total number of nodes that need to be evaluated. In fact, (RMP)$_n$ can be solved at any column generation iteration, but this is typically not worth the computational effort.
\end{remark}

\begin{remark}
For ease of exposition, a global set of columns, $\widehat{\mathcal{K}}$, is used in Algorithm \ref{alg:BnP}. However, it is usually beneficial to consider node-specific column sets. For example, branching can render some of the columns infeasible, in which case they should be removed from the set.
\end{remark}

\section{Decomposable Pricing Problems}
\label{sec:Decompose}

Most problems of interest that can be effectively solved using branch-and-price have a decomposable structure. Specifically, this means that the pricing problem decomposes into multiple smaller subproblems that can be solved independently and in parallel. Such an MINLP is of the following form:
\begin{subequations}
\label{eqn:DecompP}
\begin{align}
\minimize_{x,y,z} \quad & c^{\top} x + \sum_{i \in \mathcal{I}} f_i(y_i,z_i) \\
\st \quad & A x + \sum_{i \in \mathcal{I}} D_i \, y_i \geq b \label{eqn:DecompComplicating} \\
& g_i(y_i,z_i) \leq 0 \quad \forall \, i \in \mathcal{I} \\
& y^{\min}_i \leq y_i \leq y^{\max}_i \quad \forall \, i \in \mathcal{I} \\
& x \in \mathbb{R}_+^m \times \mathbb{Z}_+^{\bar{m}} \\
& y_i \in \mathbb{Z}^{p_i}, \, z_i \in \mathbb{R}^{q_i} \times \mathbb{Z}^{\bar{q}_i} \quad \forall \, i \in \mathcal{I},
\end{align}
\end{subequations}
where $\mathcal{I}$ denotes the set of subproblems. One can see that the problem decomposes into $|\mathcal{I}|$ independent subproblems if constraints \eqref{eqn:DecompComplicating} are removed.

%

Following an analogous derivation as in Section \ref{sec:Reformulation}, we reformulate \eqref{eqn:DecompP} into the following master problem:
\begin{subequations}
\begin{align}
v^{\mathrm{MP}} = \min_{x,\lambda} \;\; & c^{\top} x + \sum_{i \in \mathcal{I}} \sum_{k \in \mathcal{K}_i} \lambda_{ik} \, \bar{f}_{ik} \\
\mathrm{s.t.} \;\; & A x + \sum_{i \in \mathcal{I}} D_i \sum_{k \in \mathcal{K}_i} \lambda_{ik} \, \bar{y}_{ik} \geq b \\
& \sum_{k \in \mathcal{K}_i} \lambda_{ik} = 1 \quad \forall \, i \in \mathcal{I} \\
& x \in \mathbb{R}_+^m \times \mathbb{Z}_+^{\bar{m}} \\
& \lambda_{ik} \in \{0,1\} \quad \forall \, i \in \mathcal{I}, \, k \in \mathcal{K}_i,
\end{align}
\end{subequations}
where $\mathcal{K}_i$ denotes the set of feasible columns in subproblem $i$, and each $k \in \mathcal{K}_i$ is associated with a solution $\bar{y}_{ik}$ and an optimal cost $\bar{f}_{ik}$. The corresponding pricing subproblem $i$ is as follows:
\begin{subequations}
\begin{align}
\zeta_i := \min_{y_i,z_i} \;\; & f_i(y_i,z_i) - \pi^{\top} D_i \, y_i - \mu_i \\
\mathrm{s.t.} \;\; & g_i(y_i,z_i) \leq 0 \\
&y^{\min}_i \leq y_i \leq y^{\max}_i \\
& y_i \in \mathbb{Z}^{p_i}, \, z_i \in \mathbb{R}^{q_i} \times \mathbb{Z}^{\bar{q}_i}.
\end{align}
\end{subequations}

In the decomposable case, the pricing subproblems are solved independently and columns are generated for each subproblem. The optimal value $v^{\overline{\mathrm{MP}}}$ can then be bounded as follows:
\begin{equation}
v^{\overline{\mathrm{RMP}}} + \sum_{i \in \mathcal{I}} \zeta_i \leq v^{\overline{\mathrm{MP}}} \leq v^{\overline{\mathrm{RMP}}}.
\end{equation}
Note that not all $\zeta_i$ have to be nonpositive, but their sum will be.

\section{Computational Experiments}
\label{sec:Experiments}

In this section, we present computational results for four different case studies in which multiple model instances are solved. We compare the computational performance between solving the full-space model using the off-the-shelf global MINLP solver BARON 19.7.9 \cite{Kilinc2018} and solving the problem using the proposed branch-and-price algorithm. All instances are solved using 25 cores on the Mesabi cluster of the Minnesota Supercomputing Institute, a Linux cluster equipped with a set of 2.5 GHz Intel Hasewell E5-2680v3 processors. For solving the full-space models without decomposition, BARON is used with the threads option set to 25, allowing for the MILP subsolver (CPLEX) to make use of parallelization. When problems are solved using branch-and-price, CPLEX 17.1 is used to solve the restricted master problems, while BARON with a single thread is used to solve the pricing subproblems in parallel. All models were implemented in Julia using the JuMP modeling language \cite{lubin15}. The stopping criteria used for all instances are a 0.1\% optimality gap and a time limit of $10^4$ seconds. Model formulations for the full-space problem, master problem, and pricing subproblems for all case studies presented are provided in the supplementary material. 

\subsection{Example 1: Cutting Circles from Rectangles}

This problem considers cutting out a set of circles $\mathcal{I}$ of different sizes from a set of given rectangles $\mathcal{R}$, which are also of different sizes. It decides which rectangles should be used and which circles should be cut from these rectangles in order to minimize the total trim loss, which is the area of the used rectangles that was not cut into circles. More generally, it is a simple nonlinear example of a generalized assignment problem, where tasks need to be assigned to different available resources, and can also be thought of as a nonlinear example of the cutting stock problem, one of the motivating examples for the development of branch-and-price \cite{gilmore61}. This problem was originally considered in \cite{Kallrath2009}, where it is solved using column enumeration, which differs from the column generation approach used in this paper in that the master problem \eqref{eqn:MP} is solved considering the full set of columns $\mathcal{K}$, instead of a dynamically updated subset of columns.

This problem decomposes into one pricing subproblem for each rectangle $r\in\mathcal{R}$, which decides which circles should be cut from each respective rectangle. These decisions are returned as a column to the master problem, which decides which rectangles should actually be used given complicating assignment constraints that state that each circle must be assigned to exactly one rectangle. For this case study, the master problem's convexity constraint \eqref{eqn:ConvexityCon} can be relaxed to an inequality, where it is possible to not choose any column for a given rectangle if no circles are cut from that rectangle in the optimal solution. The overall size of the problem can be increased in two ways: increasing the number of rectangles, which increases the number of subproblems in the branch-and-price algorithm and should have minimal effect on solution time if the solution of the subproblems is fully parallelized, and increasing the number of circles, which increases the size of subproblems in the branch-and-price algorithm and should result in an increase in solution time. Computational results for 25 instances of varying problem size are presented in Table \ref{tab:Ex1Results}.

\begin{table}
	\caption{Computational results for 25 instances of Example 1. Optimality gap after 10,000\,s is reported if the problem is not solved within this time limit.}
	\label{tab:Ex1Results}
	\scriptsize
	\begin{tabular}{|c|c||c|c|c|c|c||c|c|}
		\hline&&Branch-&Column&Number&Branch-&Branch-&Full-&Full-\\
		$|\mathcal{I}|$&$|\mathcal{R}|$&and-Price&Generation&of Columns&and-Price&and-Price&Space&Space\\
		&&Nodes&Iterations&Generated&Time (s) or \textit{*Gap* (\%)}&Objective& Gap (\%)& Objective\\
		\hline6&20&3&7&122&335&10.2&95.1&11.2\\
		\hline6&40&1&5&185&377&9.4&98.5&12.4\\
		\hline6&60&1&5&307&576&10.3&99.9&13.0\\
		\hline6&80&3&7&448&946&9.5&N/A&not found\\
		\hline6&100&1&5&507&733&8.0&99.5&12.5\\
		\hline7&20&9&19&161&3229&11.9&98.9&14.3\\
		\hline7&40&13&21&323&4047&10.2&99.3&12.1\\
		\hline7&60&1&4&365&1489&9.5&N/A&not found\\
		\hline7&80&3&6&509&2391&9.9&N/A&not found\\
		\hline7&100&1&3&674&1112&9.3&N/A&not found\\
		\hline8&20&1&5&182&1390&10.5&98.7&13.2\\
		\hline8&40&1&4&309&1544&10.2&99.5&13.0\\
		\hline8&60&1&6&477&2562&9.9&N/A&not found\\
		\hline8&80&9&12&639&\textit{*0.67*}&9.2&N/A&not found\\
		\hline8&100&7&8&768&\textit{*2.3*}&9.3&N/A&not found\\
		\hline9&20&7&20&281&9962&9.9&99.6&19.2\\
		\hline9&40&11&20&421&\textit{*4.7*}&9.8&N/A&not found\\
		\hline9&60&5&10&595&\textit{*3.3*}&10.2&N/A&not found\\
		\hline9&80&1&6&715&\textit{*48.5*}&9.6&N/A&not found\\
		\hline9&100&5&7&986&\textit{*1.7*}&9.2&N/A&not found\\
		\hline10&20&1&8&254&2075&9.7&N/A&not found\\
		\hline10&40&1&9&459&9878&8.6&99.9&16.0\\
		\hline10&60&3&8&698&\textit{*6.0*}&9.1&N/A&not found\\
		\hline10&80&1&6&859&\textit{*64.8*}&8.9&N/A&not found\\
		\hline10&100&1&8&1056&\textit{*78.8*}&8.7&N/A&not found\\
		\hline
	\end{tabular}
\end{table}

For this case study, a feasible initial set of columns is required to ensure feasibility of the master problem in early iterations. The initial columns chosen are feasible (but likely suboptimal) cases where a single circle is cut from a rectangle, as optimal column costs can be trivially calculated for these cases as the rectangle area minus the circle area. The ``number of columns generated" information presented in Table \ref{tab:Ex1Results} includes these initial columns.
When solving the full-space problem, for many of the instances with a large number of rectangles, BARON is unable to determine an upper bound on the objective within the time allotted. Overall, the performance of the branch-and-price algorithm is clearly superior to solving the full-space problem in BARON for the instances tested in this case study, as the full-space model cannot be solved to optimality in under 10,000\,s for any instance, and only returns solutions with very large optimality gaps after that time. In contrast, 16 of the 25 instances are solved to global optimality using branch-and-price. As expected, the solution times for the branch-and-price algorithm increase notably with the number of circles, or subproblem size, but less so with the number of rectangles, or number of subproblems. Note that since we only utilize 25 processors, the branch-and-price algorithm is not fully parallelized for instances with more than 25 subproblems. Also note that 12 of the 25 instances tested in this case study require branching to find the global optimum; here, the branching strategy used is to find the circle of the largest size with a non-integer assignment to a rectangle, and to branch on this variable. 

\subsection{Example 2: Multiperiod Capacity Planning with Congestion Effects}

This problem considers capacity planning at a facility that can use a set of candidate production units $\mathcal{I}$ to produce a set of demanded products $\mathcal{J}$. Demands for these products change over time throughout a set of time periods $\mathcal{T}$. Demands are placed into a production queue, and congestion within this queue results in uncertainty in lead times for meeting demands. It decides which units to build during which time periods in order to meet demands at minimum cost. Nonlinear chance constraints are used to ensure lead times remain below a given threshold with some probability. This problem is a multiperiod extension of the capacity planning with congestion problem considered in \cite{Rajagopalan2001}.

In extending the original formulation to the multiperiod case, we define two different sets of variables which are related to a unit being available for use: $z_{it}$, which corresponds to the number of units $i$ newly built during time period $t$, and $y_{it}$, which corresponds to the number of units $i$ that are actively used during time period $t$. These variables are related according to the following conservation equations:
\begin{equation}
\label{eqn:unitcons}
y_{it}\leq y_i^0+\sum_{t'=1}^tz_{it'}\quad\forall\,i\in\mathcal{I},\,t\in\mathcal{T},
\end{equation}
where $y_{i}^0$ is the number of units $i$ initially present. This problem decomposes into one pricing subproblem for each time period $t\in\mathcal{T}$, which decides the assignment of orders to production units such that demands are met within the respective time period. Subproblems are linked to the master problem via the number of units actively used in each time period, and the master problem then determines which units should be built when subject to constraints \eqref{eqn:unitcons}. We again increase the overall size of the problem in two ways: by increasing the number of subproblems (increasing $|\mathcal{T}|$) and by increasing the size of subproblems (increasing $|\mathcal{I}|$ and $|\mathcal{J}|$). Computational results for 20 instances of varying problem size are presented in Table \ref{tab:Ex2Results}.

\begin{table}
	\caption{Computational results for 25 instances of Example 2. Optimality gap after 10,000\,s is reported if the problem is not solved within this time limit.}
	\label{tab:Ex2Results}
	\scriptsize
	\begin{tabular}{|c|c|c||c|c|c|c|}
		\hline&&&Column&Number&Branch-&Branch-\\
		$|\mathcal{I}|$&$|\mathcal{J}|$&$|\mathcal{T}|$&Generation&of Columns&and-Price&and-Price\\
		&&&Iterations&Generated&Time (s) or \textit{*Gap* (\%)}&Objective ($*10^3$)\\
		\hline5&10&6&10&50&338&5.8\\
		\hline5&10&12&19&138&917&5.8\\
		\hline5&10&24&22&196&1741&9.0\\
		\hline5&10&36&28&233&2682&10.9\\
		\hline8&15&6&34&165&5886&7.7\\
		\hline8&15&12&27&303&\textit{*0.36*}&9.6\\
		\hline8&15&24&18&346&\textit{*18.8*}&14.0\\
		\hline8&15&36&10&313&\textit{*46.4*}&15.4\\
		\hline10&20&6&2&18&389&22.7\\
		\hline10&20&12&8&57&1970&19.6\\
		\hline10&20&24&6&83&3874&26.0\\
		\hline10&20&36&7&114&3661&22.4\\
		\hline12&25&6&2&15&485&22.7\\
		\hline12&25&12&2&37&448&25.7\\
		\hline12&25&24&7&78&3482&27.8\\
		\hline12&25&36&3&106&2145&31.8\\
		\hline15&30&6&3&15&772&30.7\\
		\hline15&30&12&3&39&705&38.4\\
		\hline15&30&24&13&125&4652&29.8\\
		\hline15&30&36&3&90&3434&37.3\\
		\hline
	\end{tabular}
\end{table} 

To obtain an initial set of columns for the branch-and-price algorithm, the pricing subproblems were solved in the case where all dual variables are set to zero. We do not explicitly define starting values for the variables in either the initial subproblems or in the full-space problem. While subproblems are solved relatively quickly, we note that in all instances tested, we are unable to obtain even a feasible solution when solving the full-space problem using BARON. For this case study it is again seen that the branch-and-price approach is superior, with feasible solutions found in all instances and globally optimal solutions found in 17 out of the 20 instances tested. Interestingly, we find that solution times tend to increase more strongly with the number of subproblems than the size of subproblems, which is the opposite of what is expected. We note that this is likely due to the fact that instances with smaller subproblem sizes tend to require more column generation iterations to converge, and that the time required per iteration of column generation seems to better conform to the expected trends. We also note that in all instances tested, no branching was needed to solve the problems as the root node solutions of the relaxed master problem are integer. This observation will also be true in the next two case studies.

\subsection{Example 3: Multiscenario Synthesis of Integrated Water Networks}

This problem considers the optimal design of a waste water network with process units that add impurities to the water, as well as remediation units that remove these impurities. It considers uncertainty in the amount of impurities added by the process units and removed by remediation units, which are manifested through a set of different scenarios $\mathcal{S}$. It decides which pipes should be used to connect different units and which remediation units should be installed in order to satisfy constraints on the allowable amount of impurities in the system outputs at minimum cost. The problem considered here is an adaptation of the original problem presented in \cite{Karuppiah2008}, using the same network structure but assuming that allowable unit and pipe sizes are drawn from a discrete set $\mathcal{V}$, rather than being allowed to vary continuously. 

To solve this problem using branch-and-price, we introduce copies of the design variables for each scenario $s\in\mathcal{S}$. These variables are constrained by non-anticipativity constraints, which take into account the fact that we do not know \textit{a priori} which of the scenarios will actually be realized when making design decisions, and as such designs must be the same for all scenarios:
\begin{equation}
\label{eqn:nonant}
y_{ivs}=y_{ivs'}\quad\forall\,i\in\mathcal{I},\,v\in\mathcal{V},\,(s,s')\in\mathcal{P},
\end{equation}
where $y_{ivs}$ refers to a design decision to build unit or pipe $i$ of size $v$ in scenario $s$, and $\mathcal{P}$ denotes the set of pairs of indistinguishable scenarios. Note that the set $\mathcal{I}$ refers to all different pipes and units. The problem then decomposes into one pricing subproblem for each scenario $s\in\mathcal{S}$, which decides the optimal network design for a specific scenario. The master problem then considers the probabilities of all scenarios and finds a design which minimizes the expected cost subject to non-anticipativity constraints \eqref{eqn:nonant}. The number of subproblems is increased by increasing $|\mathcal{S}|$, and the subproblem size is increased by increasing $|\mathcal{V}|$. Computational results for 20 instances of varying size are presented in Table \ref{tab:Ex3Results}.

\begin{table}
	\caption{Computational results for 20 instances of Example 3. Optimality gap after 10,000\,s is reported if the problem is not solved within this time limit.}
	\label{tab:Ex3Results}
	\scriptsize
	\begin{tabular}{|c|c||c|c|c|c|c||c|c|}
	\hline&&Column&Number&Infeasible&Branch-&Branch-&Full&Full\\
	$|\mathcal{V}|$&$|\mathcal{S}|$&Generation&of Columns&Columns&and-Price&and-Price&Space&Space\\
	&&Iterations&Generated&Generated&Time (s) or \textit{*Gap* (\%)}&Objective ($*10^5$)& Gap (\%)& Objective($*10^5$)\\
	\hline1&5&2&3&1&172&6.9&37.2&6.9\\
	\hline1&10&7&39&1&1451&6.6&73.1&6.9\\
	\hline1&15&7&34&1&2197&6.9&N/A&not found\\
	\hline1&20&7&38&1&1537&6.6&N/A&not found\\
	\hline1&25&3&30&1&618&6.7&N/A&not found\\
	\hline2&5&12&47&1&2857&6.7&68.9&6.7\\
	\hline2&10&16&114&1&4909&6.5&70.6&6.5\\
	\hline2&15&9&66&1&5254&6.7&71.2&6.7\\
	\hline2&20&17&116&1&\textit{*1.9*}&6.7&72.5&6.7\\
	\hline2&25&14&104&1&9456&6.7&72.4&6.7\\
	\hline3&5&19&91&3&\textit{*13.6*}&6.9&70.9&7.2\\
	\hline3&10&13&125&7&\textit{*38.8*}&6.7&73.6&6.9\\
	\hline3&15&11&137&5&\textit{*46.2*}&6.9&N/A&not found\\
	\hline3&20&8&127&2&\textit{*28.1*}&6.6&N/A&not found\\
	\hline3&25&6&129&1&\textit{*34.3*}&6.7&N/A&not found\\
	\hline4&5&25&73&29&\textit{*52.0*}&6.7&74.8&7.0\\
	\hline4&10&14&99&19&\textit{*12.7*}&6.8&N/A&not found\\
	\hline4&15&10&121&8&\textit{*41.5*}&6.7&N/A&not found\\
	\hline4&20&8&122&18&\textit{*15.0*}&6.6&N/A&not found\\
	\hline4&25&12&227&22&\textit{*33.2*}&6.6&N/A&not found\\
	\hline
	\end{tabular}
\end{table}

For this problem, the knowledge of the non-anticipativity constraints is used to modify the algorithm for generating columns. Whenever a column with negative reduced cost is generated for one scenario, a cost for that column from each scenario is found by solving each pricing subproblem with design variables fixed at the value given by the new column. As such, this problem only consists of one set of columns, instead of multiple sets corresponding to each respective subproblem. Initial columns are generated using this strategy with dual variables equal to zero. This approach ensures that the master problem is always feasible. However, it is possible that columns generated from one scenario may be infeasible for another. These infeasible columns are kept in a separate set and integer cuts are added to the subproblems to ensure that they are not generated again. Note that the number of infeasible columns increases with the number of discrete unit sizes, which makes sense as this introduces a larger number columns that may be infeasible for a scenario. Again, the branch-and-price algorithm outperforms solving the full-space model using BARON in all instances. However, the time required to solve the problem clearly scales very poorly with subproblem size in this case, to the point where when $|\mathcal{V}|\geq 3$, no instance is solved to optimality within the time limit. Branch-and-price still finds a feasible solution in all cases, which is not the case when solving the full-space problem, which only finds a feasible solution in 10 out of 20 instances.

\subsection{Example 4: Multiscenario Design of Multiproduct Batch Plants}

This problem considers the optimal design of a batch plant with a set of stages $\mathcal{J}$ that can produce multiple products. It considers uncertainty in the demands for each product, as well as in the operating parameters of each of the different batch stages, which are manifested through a set of different scenarios $\mathcal{S}$. It decides the number of batch units to include at each stages, as well as the batch sizes and times for each product, such that demands are met at minimum cost. The problem presented here is an adaptation of a problem presented in \cite{Grossmann1979}, assuming that the volume of the batch unit is fixed instead of being a decision variable, and adding a nonlinear operating cost term to the objective function.

To solve this problem using branch-and-price, we introduce copies of the variable corresponding to number of batch units built for each scenario $s\in\mathcal{S}$. Like in the previous example, these variables are constrained by non-anticipativity constraints. As such, we use the same modified algorithm for generating columns. To combat the possibility of generating infeasible columns, scenarios are grouped together into a set of scenario bundles $\mathcal{B}$ such that each bundle roughly contains an equal number of ``easy-to-satisfy" scenarios where demands are low and ``hard-to-satisfy" scenarios where demands are high. Additionally, we add cuts to all subproblems when an infeasible column is generated which state that if a certain design is infeasible, any new design must build more batch units than in the infeasible design in at least one stage:
\begin{align}
&N_j+N^{\mathrm{inf}}_{jc}z_{jc}\geq N^{\mathrm{inf}}_{jc}+1\quad\forall\,j\in\mathcal{J},c\in\mathcal{C}^{\mathrm{inf}},\label{ic1}\\
&\sum_{j\in\mathcal{J}}z_{jc}\leq|\mathcal{J}|-1\quad\forall\,c\in\mathcal{C}^{\mathrm{inf}},\label{ic2}
\end{align}
where $N_j$ is the number of units built for stage $j$, $N_{jc}^{\mathrm{inf}}$ is the number of units built for stage $j$ in infeasible column $c$, $z_{jc}$ is a binary variable which is $0$ only if $N_j > N_{jc}^{\mathrm{inf}}$, and $\mathcal{C}^{\mathrm{inf}}$ is a set of infeasible columns. The problem then decomposes into one pricing subproblem for each scenario bundle $b\in\mathcal{B}$, which determines the optimal system design for the corresponding group of scenarios. The master problem then considers the probabilities of all scenarios and finds a design which minimizes the expected cost subject to non-anticipativity constraints. The number of subproblems is increased by increasing $|\mathcal{B}|$, where each bundle contains 5 scenarios, and the size of subproblems is increased by increasing $|\mathcal{J}|$. Computational results for 20 instances of varying size are presented in Table \ref{tab:Ex4Results}.

\begin{table}
		\caption{Computational results for 20 instances of Example 4. Optimality gap after 10,000\,s is reported if the problem is not solved within this time limit.}
	\label{tab:Ex4Results}
	\scriptsize
	\begin{tabular}{|c|c||c|c|c|c|c||c|c|}
		\hline&&Column&Number&Infeasible&Branch-&Branch-&Full&Full\\
		$|\mathcal{J}|$&$|\mathcal{B}|$&Generation&of Columns&Columns&and-Price&and-Price&Space&Space\\
		&&Iterations&Generated&Generated&Time (s) or \textit{*Gap* (\%)}&Objective ($*10^5$)& Gap (\%)& Objective($*10^5$)\\
		\hline3&10&9&47&0&1934&8.3&1.11&8.3\\
		\hline3&20&9&42&25&3120&8.6&8.03&8.6\\
		\hline3&30&7&58&16&2471&8.4&1.31&8.4\\
		\hline3&40&5&45&0&1349&8.5&1.44&8.5\\
		\hline3&50&5&43&39&2271&8.6&0.30&8.6\\
		\hline4&10&39&57&26&\textit{*0.13*}&16.0&3.77&16.0\\
		\hline4&20&8&41&35&2458&17.2&1.14&17.2\\
		\hline4&30&9&101&70&\textit{*1.42*}&15.6&3.12&15.6\\
		\hline4&40&12&69&91&\textit{*0.13*}&16.7&0.43&16.7\\
		\hline4&50&6&61&53&3449&16.4&0.22&16.4\\
		\hline5&10&25&119&79&\textit{*2.49*}&17.8&2.63&17.8\\
		\hline5&20&25&237&146&\textit{*1.36*}&17.5&2.99&17.5\\
		\hline5&30&15&187&96&\textit{*4.00*}&18.0&4.46&18.0\\
		\hline5&40&12&183&133&\textit{*2.47*}&17.9&1.64&17.9\\
		\hline5&50&10&139&166&\textit{*2.77*}&18.9&2.61&18.9\\
		\hline6&10&33&149&155&\textit{*5.33*}&21.9&4.03&21.9\\
		\hline6&20&32&238&301&\textit{*4.78*}&22.5&3.16&22.5\\
		\hline6&30&18&242&192&\textit{*3.16*}&23.1&3.01&23.1\\
		\hline6&40&14&212&274&\textit{*4.44*}&23.0&2.26&23.0\\
		\hline6&50&12&221&238&\textit{*3.29*}&24.4&4.86&24.3\\
		\hline
	\end{tabular}

\end{table}

For this problem, unlike the other three examples presented, the branch-and-price approach performs very similarly to solving the full-space problem using BARON, particularly as the subproblem size increases. For all problems, both approaches are able to find the same upper bounds. We note that for some instances, BARON is able to find a slightly better objective value than the branch-and-price approach, although in all instances these differences are within either the 0.1\% gap used as a global optimality stopping criterion, or within the reported gap after 10,000\,s. It is also apparent that as the number of stages increases, the number of infeasible columns generated also increases, contributing to the added difficulty for solving this problem using branch-and-price. However, we do note that the number of infeasible columns generated is reduced by bundling columns and adding infeasibility cuts, as this number can be as much as an order of magnitude greater when these options are not used.

\section{Conclusions}
\label{sec:Conclusions}

In this work, we have applied branch-and-price to a class of nonconvex MINLPs with linear complicating constraints and integer linking variables. We exploit the structure of the problem to construct a Dantzig-Wolfe reformulation via a discretization approach, which then allows the problem to be solved using a branch-and-price scheme. This approach is especially effective in cases where the pricing problem decomposes into multiple small subproblems such that solving each subproblem using a global MINLP solver is considerably more tractable than solving the original full-space MINLP. Through several case studies, we have shown that many relevant problems directly fall or can be reformulated into the given class of MINLPs, and have demonstrated the computational feasibility of the proposed algorithm. In most tested model instances, the branch-and-price algorithm clearly outperforms solving the full-space problem directly using a global MINLP solver, often achieving orders-of-magnitude speedups.

\section*{Acknowledgments}
We gratefully acknowledge financial support from the University of Minnesota and the Minnesota Supercomputing Institute (MSI) at the University of Minnesota for providing resources that contributed to the research results reported within this paper. We also thank Angela Flores-Quiroz for insightful discussions on our work.

\bibliographystyle{newapa}
\bibliography{library}

\end{document}